\documentclass{article}
\usepackage{amssymb}
\usepackage[cp1251]{inputenc}
\usepackage[english]{babel}
\usepackage{amsmath}

\def\mapr#1{\smash{\mathop{\buildrel{#1}\over\longrightarrow}}}

\def\mapu#1{\Big\uparrow\rlap{$\vcenter{\hbox{$#1$}}$}}
\def\cK{{\cal K}}
\def\Z{{\bf Z}}
\def\K{{\bf K}}
\def\Q{{\bf Q}}
\newtheorem{thr}{Theorem}
\def\({\left(}
\def\){\right)}
\def\<{\langle}
\def\>{\rangle}

\title{Construction
of Fredholm representations and a modification of the Higson-Roe corona
}
\author{
A.S.Mishchenko and N.Teleman}
\date{}
\begin{document}
\parindent 0.0cm
\maketitle

\begin{abstract}
The Fredholm representation theory is well adapted to construction
of homotopy invariants of non simply connected manifolds on the
base of generalized Hirzebruch formula
\begin{equation}\label{eq1.1}
[\sigma (M)]=\langle L(M)\hbox{ch}_{A}f^{*}\xi ,[M]\rangle \in \mbox{\large\bf
K}^{0}_{A}(\hbox{pt})\otimes \mbox{\large\bf Q},
 \end{equation}
where  $A=C^{*}[\pi ]$  is the group $C^{*}$--algebra of the group $\pi $, $\pi
=\pi _{1}(M)$. The bundle $\xi\in \mbox{\large\bf K}^{0}_{A}(B\pi )$ is
canonical $A$--bundle, generated by the natural representation $\pi \mapr{}A $
.

In \cite{Mi129-e} a natural family of the Fredholm representations was constructed
that lead to a symmetric vector bundle on completion of the fundamental group
with a modification of the Higson-Roe corona when the completion is a closed manifold.

 Here we will discuss a homology version of symmetry in the case
 when completion with a modification of the Higson-Roe corona is
 a manifold with boundary. The results were developed 
 during the visit of the first author  in Ancona on March, 2007.
 The second version is
 supplemented by details of consideration the case
 of manifolds with boundary.
\end{abstract}

The Fredholm representation theory is well adapted to construction
of homotopy invariants of non simply connected manifolds on the
base of generalized Hirzebruch formula
\begin{equation}\label{eq1.1a}
[\sigma (M)]=\langle L(M)\hbox{ch}_{A}f^{*}\xi ,[M]\rangle
\in \mbox{\large\bf K}^{0}_{A}(\hbox{pt})\otimes \mbox{\large\bf Q}
 \end{equation}
where  $A=C^{*}[\pi ]$  is the group $C^{*}$--algebra of the group
$\pi $, $\pi =\pi _{1}(M)$. The bundle $\xi\in \mbox{\large\bf
K}^{0}_{A}(B\pi )$ is canonical $A$--bundle, generated by the
natural representation $\pi \mapr{}A $ . The map $f:M \mapr{}B\pi $
induces the isomorphism of fundamental groups. The element $[\sigma
(M)]\in \mbox{\large\bf K}^{0}_{A}(\hbox{pt}) $
is generated by noncommutative signature of the manifold $M$
under exchange of rings $\mbox{\large\bf Z}[\frac{1}{2}][\pi
]\subset A$.

Let $\rho =(T_{1}, F, T_{2})$ be a Fredholm representation of the
group $\pi $, that is a pair of unitary representations
$T_{1},T_{2}:\pi \mapr{}B(H) $ and a Fredholm operator
$F:H\mapr{}H$, such that

 \begin{equation}\label{eq1.2}
FT_{1}(g)-T_{2}(g)F\in \hbox{Comp}(H),\; g\in \pi .
 \end{equation}

Changing the algebra $B(H)$ for the Calkin algebra
 ${\cal K}=B(H)/\hbox{Comp}(H)$,
one comes to the representation $\widehat\rho $ of the group
$\pi \times {\bf Z}$ to the Calkin algebra:
\begin{equation}\label{eq1.2a}
\widehat\rho:\pi \times {\bf Z}\mapr{}\cK,
\end{equation}
\begin{equation}\label{eq1.2b}
  \widehat\rho(g,n) = T_{2}(g)F^{n}= F^{n}T_{1}(g), \quad g\in\pi, \quad n\in \Z.
\end{equation}


 \begin{equation}\label{eq1.3a}
\rho_{*} :K_{A}(X) \mapr{\hbox{Id}\otimes \beta } K_{A\widehat\otimes
C(S^{1})}(X\times S^{1}) \mapr{\widehat{\rho} }K_{{\cal K}}(X\times S_{1}).
 \end{equation}
Here $\beta \in K_{C(S^{1})}(S^{1})$ is the canonical element
generated by regular representation of the group ${\bf Z}$.

Applying (\ref{eq1.3a}) to the Hirzebruch formula (\ref{eq1.1a}) one has homotopy
invariance of corresponding higher signature.

\section{Construction of Fredholm representation}\label{sec.1}


Let $T$ be the sum of finite copies of regular
representation of the group  $\pi $,  $\Phi $ be the block diagonal
operator that is defined as matrix valued function $F(g),\; g\in \pi $:
 \begin{equation}\label{eq1.3.1}
F(g):V \mapr{}V.
 \end{equation}

Let
 \begin{equation}\label{eq1.3.2}
H=\bigoplus_{g\in\pi } V_{g}, \; V_{g}\equiv V,
 \end{equation}

 \begin{equation}\label{eq1.3.3}
T_{h}: H \mapr{}H,\; V_{g}\mapr{}V_{hg}.
 \end{equation}
The condition that $\Phi $ is Fredholm operator means that

 \begin{equation}\label{eq1.4}
\|F(g)\|\leq C,\; \|F^{-1}(g)\|\leq C
 \end{equation}
 for all  $g\in \pi $ except a finite subset.
The condition (\ref{eq1.2}) means  that

 \begin{equation}\label{eq1.5}
\lim_{|g|\mapr{}\infty }\|F(g)-F(hg)\|=0.
 \end{equation}
So if the pair
 \begin{equation}\label{eq1.6.1}
\rho =(T,\Phi )
 \end{equation}
 satisfies conditions (\ref{eq1.4}), (\ref{eq1.5}), then
$\rho $ is Fredholm representation of the group
 $\pi $.
 \vskip 1cm

 Consider universal covering $\widetilde{B\pi} $
of classifying space $B\pi $ endowed with left action of the group
 $\pi $.
In correspondence to the construction by    \cite{Mi17-70e}
the vector bundle generated by the representation $\rho$ on the space
$B\pi $ can be represented as an equivariant
continuous family of Fredholm operators on the space
$E\pi =\widetilde{B\pi} $.
The property of equivariance corresponds to diagonal action on Cartesian product

 \begin{equation}\label{eq1.6.2}
T_{h}:E\pi \times H,\; (x,\xi ) \mapr{}(hx,T_{h}(\xi )).
 \end{equation}
Namely, let the space$B\pi $ be endowed with a structure
of simplicial space and  $E\pi =\widetilde{B\pi} $ be endowed with the structure
of simplicial structure derived from the covering

 \begin{equation}\label{eq1.6}
E\pi =\widetilde{B\pi} \mapr{p}  B\pi
\end{equation}

Let  $\{x_{i}\}$ be the family of vertices of
$E\pi =\widetilde{B\pi}$ , one from each of orbits
of the action of $\pi $. Then each simplex
$\sigma $ of $E\pi =\widetilde{B\pi}$  is defined completely by their vertices
\begin{equation}\label{eq1.7}
\sigma =(h_{0}x_{i_{0}},\dots, h_{n}x_{i_{n}}), \; h_{0},\dots,h_{n}\in\pi .
 \end{equation}
Any point  $x\in\sigma $ is uniquely defined as a convex
linear combination of vertices
 \begin{equation}\label{eq1.8}
x=\sum^{n}_{k=0}\lambda _{k}h_{k}x_{i_{k}}
 \end{equation}
Then the equivariant family
of Fredholm operators which corresponds to the Fredholm representation
$\rho$ (\ref{eq1.6.1}) one can define by the formula

 \begin{eqnarray}\label{eq1.9}
\lefteqn{\Phi _{x}=
\Phi _{x}(\rho)=
}\nonumber \\
& = &
\sum^{n}_{k=0}\lambda _{k}\Phi _{h_{k}x_{i_{k}}} =
\sum^{n}_{k=0}\lambda _{k}T_{h_{k}}\Phi _{x_{i_{k}}}T^{-1}_{h_{k}} =
\nonumber \\
& = &
\sum^{n}_{k=0}\lambda _{k}T_{h_{k}}\Phi T^{-1}_{h_{k}} .
 \end{eqnarray}
Hence
 \begin{equation}\label{eq1.10}
\left( \Phi _{x}\right) _{g}=
\sum^{n}_{k=0}\lambda _{k}F_{h_{k}^{-1}g} .
 \end{equation}
It is clear that the family (\ref{eq1.9}) is equivariant. Indeed,
 \begin{equation}\label{eq1.11}
hx=\sum^{n}_{k=0}\lambda _{k}hh_{k}x_{i_{k}} .
 \end{equation}
Hence
 \begin{equation}\label{eq1.12}
\Phi _{hx}=
\sum^{n}_{k=0}\lambda _{k}T_{hh_{k}}\Phi T^{-1}_{hh_{k}} =
T_{h}\left( \sum^{n}_{k=0}\lambda _{k}T_{h_{k}}\Phi T^{-1}_{h_{k}}\right) T^{-1}_{h} =
T_{h}\Phi _{x}T^{-1}_{h}.
 \end{equation}
Also it is clear that the operators
 (\ref{eq1.9}) are Fredholm by
(\ref{eq1.10}) , (\ref{eq1.5}) and (\ref{eq1.4}).

\vskip 1cm
On the other side the operators
 (\ref{eq1.3.1}) generate the continuous family

 \begin{equation}\label{eq1.13}
F_{x}: V \mapr{}V, \; x\in E\pi
 \end{equation}
using formula
 \begin{equation}\label{eq1.14}
F_{x}=\sum^{n}_{k=0}\lambda _{k}F(h^{-1}_{k}).
 \end{equation}
This family one can consider as a linear mapping of the trivial bundle:
\begin{equation}\label{eq1.15}
F_{x}:E\pi \times V \mapr{} E\pi \times V.
 \end{equation}
Consider the universal covering
 \begin{equation}\label{eq1.16}
p:E\pi \mapr{}B\pi .
 \end{equation}
Denote
 \begin{equation}\label{eq1.17}
{\cal K}^{i}(E\pi )=\lim_{\leftarrow}K_{c}^i(p^{-1}(X)),
 \end{equation}
 where the inverse limit takes by the family of
 all compact subsets $X\subset B\pi $.

 \begin{thr}\label{th.1}
The map (\ref{eq1.15}) defines the element

 \begin{equation}\label{eq1.18}
F(\rho )\in {\cal K}^{0}(E\pi ).
 \end{equation}
 \end{thr}
\vskip 1cm

Consider the direct image of the bundle
(\ref{eq1.15}) over $B\pi $:

 \begin{equation}\label{eq1.16}
A \mapr{}B\pi ,
 \end{equation}
where the fibre is the direct sum of the fibers of the bundle (\ref{eq1.15})
over each orbit of the action of the group $\pi $ in the space  $E\pi $.
The total space $A$ is defined as
 \begin{equation}\label{eq1.17}
A=\{(u,\xi ):u\in B\pi , \xi \in\bigoplus_{x\in u}(x\times V)\}.
 \end{equation}
 Let
 \begin{equation}\label{eq1.18}
\widetilde A \mapr{}  E\pi
 \end{equation}
 be the inverse image of the bundle (\ref{eq1.16}).
The total space
$\widetilde A$ is defined as
 \begin{equation}\label{eq1.19}
\widetilde A=\{(x,\xi ):x\in E\pi , \xi \in \bigoplus_{y\in[x]}(y\times V)\}=
\{(x,\xi ),x\in E\pi , \xi \in \bigoplus_{g\in\pi }(gx\times V)\}.
 \end{equation}
Define the action of the group $\pi $ on the total space
$\widetilde A$ by the formula
 \begin{eqnarray*}\label{eq1.20}
f_{h}(x,\xi )
& = & (hx,\eta ),
\nonumber \\
\xi
& = &
\oplus \xi _{g} \in \bigoplus_{g\in\pi }(gx\times V),
\nonumber \\
\eta
& = &
\oplus\eta_{g}  \in \bigoplus_{g\in\pi }(ghx\times V),
\nonumber \\
\eta _{g}
& = &
\xi _{gh}.
\end{eqnarray*}
It is clear that
 \begin{equation}\label{eq1.21}
A=\widetilde A/\pi .
 \end{equation}
On the other side there is an isomorphism
between the bundle (\ref{eq1.18}) and the bundle (\ref{eq1.6.2}):

 \begin{equation}\label{eq1.22}
\varphi : E\pi \times \bigoplus_{g\in\pi } V_{g} \mapr{}\widetilde A,
 \end{equation}
 \begin{equation}\label{eq1.23}
\varphi (x, \oplus \xi _{g})=(x,\oplus \xi _{g^{-1}}.
 \end{equation}
This isomorphism is equivariant. The map (\ref{eq1.15}) goes to the map
of the direct image as the mapping
 \begin{equation}\label{eq1.24}
\widetilde F:\widetilde A \mapr{}\widetilde A,
 \end{equation}

 \begin{eqnarray*}\label{eq1.25}
\widetilde F (x, \oplus \xi _{g}) & = & \left( x, \oplus F_{gx}(\xi
_{g})\right)  = \nonumber \\ & = & \left( x, \oplus \sum^{n}_{k=0}\lambda
_{k}F_{h^{-1}_{k}g^{-1}} (\xi _{g})\right) .
 \end{eqnarray*}
 \begin{equation}\label{eq1.25}
\widetilde F (x, \oplus \xi _{g})=(x, \oplus F_{gx}(\xi _{g})
 \end{equation}

It is clear that the map (\ref{eq1.24}) goes to (\ref{eq1.9})
under the isomorphism (\ref{eq1.22}).

So the following theorem holds:
 \begin{thr}\label{th.2}
 Consider the Fredholm representation of the group $\pi $
 of the form (\ref{eq1.6.1}). Let
$\xi _{\rho }\in \mbox{\large\bf K}(B\pi )
$
 be the element defined by the mapping  (\ref{eq1.9}) .
Then
 \begin{equation}\label{eq1.27}
p_{!}(F(\rho )) =\xi _{\rho }\in \mbox{\large\bf K}^{0}(B\pi )   ,
 \end{equation}
where
 \begin{equation}\label{eq1.28}
p_{!}: {\cal K}^{0}(E\pi ) \mapr{} \mbox{\large\bf K}^{0}(B\pi )
 \end{equation}
 is the direct image in  $K$--theory.
 \end{thr}
 \vskip 1cm

Consider the action of the group $\pi $ on the Cartesian product
$E\pi \times V$ as the left action on the first factor and identical
on the second one.

Consider on the space ${E\pi} $ a metric with the property
 \begin{equation}\label{eq1.6}
 r(xg,yg)\mapr{}0,\; |g|\mapr{}\infty.
 \end{equation}
Let  $\overline{E\pi }$ be the completion of the space
 ${E\pi} $ (with respect to the metric $r$).
Then any continuous mapping
 \begin{equation}\label{eq1.7}
f:(\overline{E\pi },\;\overline{E\pi }\backslash{E\pi})\mapr{}
(B(V),\;U(V))
 \end{equation}
defines the continuous family
of the Fredholm representations  $\rho(x)$, $x\in E\pi$.

By the theorem \ref{th.1} the family $\rho(x)$ generates the equivariant
family
\begin{equation}\label{eq1.7a}
F_{x,y}:E\pi \times E\pi\times V \mapr{} E\pi \times E\pi \times V.
\end{equation}
and therefore the element

 \begin{equation}\label{eq1.7b}
F(\rho(x) )\in {\cal K}^{0}\((E\pi \times E\pi)/\pi\).
 \end{equation}

Let
 \begin{equation}\label{eq1.7c}
p_{!}: {\cal K}^{0}\((E\pi\times E\pi)/\pi\) \mapr{}
\mbox{\large\bf K}^{0}(B\pi\times B\pi )
 \end{equation}
be the direct image in  $K$--theory. Then
\begin{equation}\label{eq1.7d}
p_{!}(F(\rho(x) )) =\xi _{\rho(x)}\in \mbox{\large\bf K}^{0}().
\end{equation}

The symmetric property of the element $\xi _{\rho(x)}$ holds:
\begin{equation}\label{eq1.7e}
  (1\otimes u)\xi _{\rho(x)}=
  (u\otimes 1)\xi _{\rho(x)}\in \K^{0}(B\pi\times B\pi),\quad u\in \K^{0}(B\pi).
\end{equation}

\section{Symmetric cohomology classes in $H^*(M\times M)$}
In the case when the space $B\pi$ is a compact manifold
and the space  ${E\pi}$ is compactified to the disk with extension
of the action of $\pi$, we obtain new proof of the Novikov conjecture in the case
 \cite{FaHs-81e}.

For that consider a closed orientable compact manifold $M$ and a cohomology class
$w\in H^*(M\times M; \Q)$. Assume that $w$ satisfies a symmetric condition:

\begin{equation}\label{s1}
 w\cdot (1\otimes x) = (x\otimes 1)\cdot w, \quad x\in H^*(M;\Q).
\end{equation}

Our aim is to describe such symmetric elements $w$. Let $x_{i}$, $0\leq i\leq
N$ be a (homogenious) basis in $H^*(M;\Q)$, $x_0=1\in H^{0}(M;\Q)$, $x_N\in
H^n(M;\Q)$, $\dim M =n$, $\<x_N,[M]\>=1$.

Then the multiplication tensor $\lambda^k_{ij}$ is defined by the formula
\begin{equation}
    x_i\cdot x_j = \lambda^k_{ij}x_k,
\end{equation}
\begin{equation}
    \lambda^k_{i0}=\lambda^k_{0i}= \delta^k_i,
\end{equation}
\begin{equation}
    \lambda^N_{ij}=\<x_i\cdot x_j,[M]\>.
\end{equation}

Associativity of the multiplication means that
\begin{equation}
    (x_i\cdot x_j)\cdot x_k =x_i\cdot (x_j\cdot x_k),
\end{equation}
that is
\begin{equation}
\begin{array}{c}
    \lambda^l_{ij}\lambda^s_{lk}x_s=(\lambda^l_{ij}x_l)\cdot x_k=\phantom{aaaaaaaaaaaaaaaa}\\\\
        =(x_i\cdot x_j)\cdot x_k =x_i\cdot (x_j\cdot x_k)=\\\\
    \phantom{aaaaaaaaaaaaaaaa}  =x_i\cdot (\lambda^l_{jk}x_l)= \lambda^s_{il}\lambda^l_{jk}x_s,
\end{array}
\end{equation}
that is

\begin{equation}\label{21}
    \lambda^l_{ij}\lambda^s_{lk}=\lambda^s_{il}\lambda^l_{jk}.
\end{equation}

Represent the element $w$ in the form
\begin{equation}
    w = \mu^{ij}x_i\otimes x_j.
\end{equation}

Then the condition (\ref{s1}) can be written as
\begin{equation}
    \mu^{il}x_i\otimes x_l\cdot x_k = \mu^{lj}x_k\cdot x_l\otimes x_j
\end{equation}

or
\begin{equation}
    \mu^{il}x_i\otimes (\lambda^j_{lk}x_j) =
    \mu^{lj}(\lambda^i_{kl}x_i)\otimes x_j,
\end{equation}

or
\begin{equation}\label{s2}
    \mu^{il}\lambda^j_{lk} =
    \mu^{lj}\lambda^i_{kl}.
\end{equation}

Assume that we have the case
\begin{equation}\label{56}
    \mu^{Nj}=\mu^{jN}=\delta^j_0.
\end{equation}

Then from (\ref{s2}) one has
\begin{equation}\label{s3}
    \mu^{il}\lambda^N_{lk} =
    \mu^{lN}\lambda^i_{kl}.
\end{equation}

or

\begin{equation}\label{s4}
    \mu^{il}\lambda^N_{lk} =
    \delta^l_0\lambda^i_{kl}=\lambda^i_{k0} = \delta^i_k.
\end{equation}

This means that the matrix $\|\mu^{ij}\|$ is the inverse matrix of the matrix $\|\lambda^N_{ij}\|$:
\begin{equation}\label{29}
    \|\mu^{ij}\|=\|\lambda^N_{ij}\|^{-1}.
\end{equation}

The rest of relations from (\ref{s2}) are the consequence from associativity (\ref{21}):
\begin{equation}\label{30}
\begin{array}{c}
    \lambda^N_{i'i}\mu^{il}\lambda^j_{lk} =\lambda^N_{i'i}\mu^{lj}\lambda^i_{kl},\\\\
    \delta^{i'}_l\lambda^j_{lk} =\lambda^N_{i'i}\mu^{lj}\lambda^i_{kl}, \\\\
    \lambda^j_{i'k} =\mu^{lj}\lambda^i_{kl}\lambda^N_{i'i},\\\\
    \lambda^N_{jj'}\lambda^j_{i'k} =\lambda^N_{jj'}\mu^{lj}\lambda^i_{kl}\lambda^N_{i'i},\\\\
    \lambda^N_{jj'}\lambda^j_{i'k} =\delta^{l}_{j'}\lambda^i_{kl}\lambda^N_{i'i},\\\\
    \lambda^N_{jj'}\lambda^j_{i'k} =\lambda^i_{kj'}\lambda^N_{i'i},\\\\
    \lambda^j_{i'k}\lambda^N_{jj'} =\lambda^N_{i'i}\lambda^i_{kj'},\\\\
\end{array}
\end{equation}
Compare with (\ref{21}):
\begin{equation}\label{21a}
    \lambda^l_{ij}\lambda^s_{lk}=\lambda^s_{il}\lambda^l_{jk}.
\end{equation}

As a consequence from (\ref{29}) one can obtain relations for symmetric elements of the form
\begin{equation}\label{31}
    w= (x\otimes 1)(\mu^{ij}x_i\otimes x_j)(1\otimes y)=(\mu^{ij}x_i\otimes x_j)(1\otimes xy).
\end{equation}

\section{Manifolds with boundary }

Assume now that  a closed orientable compact manifold $M$ has nonempty boundary $\partial M$.
Then one has the Poincare duality as a commutative diagram
\begin{equation}
    \begin{array}{ccccccccccc}
\cdots\mapr{}&H_{k+1}(M) &\mapr{j}& H_{k+1}(M,\partial M) &\mapr{\delta}& H_{k}(\partial M) &\mapr{i}\\\\
&\mapu{D}&&\mapu{D}&&\mapu{D} \\\\
\cdots\mapr{}&H^{n-k}(M,\partial M) &\mapr{j^*}& H^{n-k}(M) &\mapr{i^*}& H^{n-k}(\partial M) &\mapr{\delta^*}\\\\\\\\
\mapr{i}&    H_{k}(M) &\mapr{j}& H_{k}(M,\partial M)&\mapr{}\cdots\\\\
&\mapu{D}&&\mapu{D}&  \\\\
\mapr{\delta^*}&H^{n+1-k}(M,\partial M) &\mapr{j^*}& H^{n+1-k}(M)&\mapr{}\cdots\\
    \end{array}
\end{equation}

The Poincare duality has the relation with multiplication in cohomology by the formula
\begin{equation}
    \<x\wedge y, [M]\> = \(x, Dy\).
\end{equation}

Here $x\in H^*(M) \quad y\in H^*(M,\partial M)$, or $y\in H^*(M) \quad x\in H^*(M,\partial M)$
and the operation $\wedge$ defines the pairing
\begin{equation}\label{d1}
    \wedge: H^i(M)\times H^j(M,\partial M) \mapr{}H^{i+j}(M,\partial M)
\end{equation}
such the pairing (\ref{d1}) generates  the module structure over the ring $H^*(M)$ with the action
on the $H^*(M,\partial M)$:
\begin{equation}\begin{array}{l}
                   y\wedge (x_1\cdot x_2) = (y\wedge x_1)\cdot x_2 =
    \pm x_1\cdot(y\wedge x_2), \\
    y\in H^*(M), \quad x_{1}, x_{2}\in  H^*(M,\partial M).\\
                \end{array}
\end{equation}

Consider Cartesian square $M\times M$. The boundary $\partial (M\times M)$ is the manifold
which is splitted into the union
\begin{equation}
  \begin{array}{l}
    \partial (M\times M) = (M\times\partial M) \cup (\partial M \times M),\\\\
    (M\times\partial M) \cap (\partial M \times M) = \partial M \times \partial M.
    \end{array}
\end{equation}

Consider a cohomology class
$w\in H^*(M\times M, \partial M \times M;\Q)$.
This cohomology can be described as a tensor product
\begin{equation}
    H^*(M\times M, \partial M \times M;\Q) \approx H^*(M,\partial M;\Q)\otimes H^*(M;\Q).
\end{equation}

Assume that $w$ satisfies a symmetric condition:

\begin{equation}\label{39}
 w\cdot (1\otimes y) = (y\otimes 1)\cdot w
 \in H^*(M\times M, \partial M \times M), \quad y\in H^*(M;\Q).
\end{equation}

The result is similar to the manifolds without boundary:
\begin{thr}
Let $w\in H^*(M\times M, \partial M \times M)$ satisfy the symmetry condition
(\ref{39}). Let $x_i\in H^*(M,\partial M)$, $y_j\in H^*(M)$ be bases,
\begin{equation}\label{40}
    w =\mu^{ij}x_i\otimes y_j.
\end{equation}
Then
\begin{equation}\label{41}
    \|\mu^{ij}\|=\|\lambda^N_{ij}\|^{-1}.
\end{equation}
where
\begin{equation}\label{29a}
    \lambda^N_{ij}= \<y_i\wedge x_j, [M,\partial M]\>.
\end{equation}
\end{thr}

To prove this consider (homogenious)  bases in the cohomology groups
$H^*(M;\Q)$ and $H^*(M,\partial M;\Q)$:
$$x_i\in H^*(M,\partial M), \quad y_j\in H^*(M).$$

 Let  $$\begin{array}{l}
   y_0=1\in H^{0}(M;\Q)\approx \Q, \\
   x_N\in
H^n(M,\partial M;\Q)\approx\Q, \quad \dim M =n,\quad \<x_N,[M,\partial M]\>=1.
 \end{array}$$
The pairing (\ref{d1}) is defined by the formula
\begin{equation}\label{d2}
y_{i}\wedge x_{j}= \lambda^{k}_{ij}x_{k}.
\end{equation}
If $y_{i}, y_{k}\in H^*(M)$ then
\begin{equation}
y_{i}\wedge y_{k}= \nu^{s}_{ik}y_{s},
\end{equation}
such that $\nu^{s}_{0k}=\nu^{s}_{k0}=\delta^{s}_{k}.$
The property (\ref{39}) can be rewrited as
\begin{equation}
\mu^{ij}x_i\otimes y_jy_{k}=\mu^{ij}y_{k}\wedge x_i\otimes y_j
\end{equation}
or
\begin{equation}
\mu^{ij}\nu^{s}_{jk}x_i\otimes y_{s}=
\mu^{ij}\lambda^{s}_{ki}x_{s}\otimes y_j.
\end{equation}
Hence
\begin{equation}
\mu^{ij}\nu^{s}_{jk}=
\mu^{ls}\lambda^{i}_{kl}.
\end{equation}

In particular when $i=N$ one has
\begin{equation}
\mu^{Nj}\nu^{s}_{jk}=
\mu^{ls}\lambda^{N}_{kl}.
\end{equation}

Assume similar to manifolds without boundary that the element
$w$ satisfies the condition (\ref{56}):

\begin{equation}\label{56}
    \mu^{Nj}=\delta^j_0.
\end{equation}
Hence
\begin{equation}
\delta^{s}_{k}=
\mu^{ls}\lambda^{N}_{kl}.
\end{equation}

The results were partially supported by the grant of
RFBR No 05-01-00923-a, the grant of the support for Scientific
Schools No NSh-619.2003.1, and the grant of the foundation
"Russian Universities" project No. RNP.2.1.1.5055


 \end{document}